\documentclass{article}
\usepackage{amsmath}
\usepackage{graphicx}
\newcommand{\eps}{\varepsilon}
\numberwithin{equation}{section}
\title{On the Origin of the Korteweg-de Vries Equation}
\author{E.~M. de Jager}
\date{\normalsize Korteweg-de Vries Institute,
University of Amsterdam,\\
P.O.~Box 94248, 1090 GE Amsterdam, The Netherlands\\
{\tt jagves@xs4all.nl}}
\begin{document}
\maketitle

%\subsection*{Abstract}
\begin{abstract}
The Korteweg-de Vries equation has a central place in a model for
waves on shallow water and it is an example of the propagation of
weakly dispersive and weakly nonlinear waves. Its history spans a
period of about sixty years, starting with experiments of Scott
Russell in 1834, followed by theoretical investigations of, among
others, Lord Rayleigh and Boussinesq in 1871 and, finally,
Korteweg and De Vries in 1895.

In this essay we compare the work of Boussinesq and Korteweg-de
Vries, stressing essential differences and some interesting
connections. Although there exist a number of articles, reviewing
the origin and birth of the Korteweg-de Vries equations,
connections and differences, not generally known, are reported.
\vskip1\baselineskip

\noindent A.M.S.~Classification: Primary 01-02, 01A55; Secondary
76-03, 76B25, 35Q53.

\end{abstract}
\newpage

\section{Introduction}

It was in the {\it ``interest of Higher Truth'' } that professor
Martin Kruskal, at the conference in commemoration of the
centennial of the publication of the Korteweg-de Vries paper in
the Philosophical Magazine \cite{ref1}, claimed that {\it ``he,
together with professor Norman Zabusky, was the person, who more
than anyone else, resuscitated  the Korteweg-de Vries equation
after its long period of, if not oblivion, at least neglect''},
\cite{ref2}. Indeed it is well-known that in the follow-up of
their 1965 paper in the Physical Review Letters \cite{ref3},
``Interactions of Solitons in a collisionless plasma and the
recurrence of initial states'', a real explosion of research on
this and related equations appeared in the journals. Many
developments in several fields of pure and applied mathematics,
physics, chemistry, biology and engineering followed. Restricting
to mathematics we mention analysis, integrability of nonlinear
systems, Lie-algebra's, differential geometry, quantum and
statistical mechanics, \cite{ref2, ref4, ref5}.

In this essay we direct our attention to the origin of the
Korteweg-de Vries equation and its birth which has been a long
process and spanned a period of about sixty years, beginning with
the experiments of Scott-Russell in 1834 \cite{ref6}, the
investigations of Boussinesq and Rayleigh around 1870 [7--11] and
finally ending with the article by Korteweg and De Vries in 1895
\cite{ref1}.

In simplified form the Korteweg-de Vries equation reads
\begin{equation}\label{eq0}
\frac {\partial u}{\partial t} - 6 u\frac {\partial u}{\partial x}
+\frac {\partial^3 u}{\partial x^3}=0
\end{equation}
and it is the result of research concerning long waves in shallow
water;  $x$ and $t$ denote position and time and $u = u (x, t)$
the wave surface.

Nowadays, it is hard to understand for mathematicians not
specialized in fluid mechanics, that a subject as this could raise
such a wide spread interest. However, this was not the case in the
nineteenth century when the study of water waves was of vital
interest for applications in naval architecture and for the
knowledge of tides and floods. Notably in England and France much
research was spent on the study of water waves of several kinds,
in England by, among others, Scott--Russell, Airy, Stokes,
McCowan, Lord Rayleigh and Lamb and in France by Lagrange,
Clapeyron, Bazin, St.~Venant and Boussinesq.

In some treatises and textbooks on  ``soliton theory''  a short
survey of the early history is presented \cite{ref12, ref13}, but
apparently it has not been the intention of the authors to dwell
extensively upon the considerations and the mathematical analysis
of those present at the cradle of the equation that became later
known as the Korteweg-de Vries equation. Nevertheless, there are
some review papers where more specific attention has been given to
investigations related to the Korteweg-de Vries equation. We
mention in particular the reviews by Bullough \cite{ref14},
Bullough and Caudry \cite{ref15}, Miles \cite{ref16} and the
recent impressive extensive article by Darrigol \cite{ref17} and a
letter in the Notices of the A.M.S. by Pego \cite{ref18}.

In these articles the work by Boussinesq on the one side and that
of Korteweg and De Vries on the other side have been discussed.
Studying these papers, the present author became aware of some
inaccuracies regarding the relevance and the significance of the
work by Korteweg and De Vries, maybe even a slight animosity over
the priority of the discovery of the equation. For example, Miles,
Darrigol and Pego suggest that Korteweg and De Vries were
presumably unaware of the work by Boussinesq. This is to be
doubted because in their article reference has been made to the
Comptes Rendus papers by Boussinesq \cite{ref7} and St.~Venant
\cite{ref19}. Besides this, the historian B.~Willink \cite{ref20}
has presented the author with a handwritten copy of De Vries,
containing an excerpt of the paper by St.~Venant \cite{ref19} and
there appears a clear reference to the {\it  ``Essai sur la
th\'eorie des eaux courantes''} \cite{ref10}, which proves that De
Vries was certainly aware of the existence of Boussinesq's
research. Pego writes {\it  ``It is not clear why Korteweg and De
Vries thought the permanence of the solitary wave still
controversial in 1895''} \cite{ref18}. This is in contrast with
the introduction of the KdV article, where it is stated {\it
 ``They (Lord Rayleigh and McCowan) are as it seems to us, inclined
to the opinion that the wave is only stationary to a certain
approximation. It is the desire to settle this question
definitively which has led us into somewhat tedious calculations,
which are to be found at the end of our paper'' }, \cite{ref1}.

It is evident that Korteweg and De Vries, wanting to check the
theory of long waves in shallow water, use an independent
approach. It is our intention to illustrate this in the next
sections, pointing out not only differences but also close
connections in both theories. Here we already give some examples.
Boussinesq used a fixed coordinate system and Korteweg and De
Vries a coordinate system moving with the wave. The central
equations in Boussinesq's analysis are the continuity equation and
an expression for the wave velocity \cite{ref9}, whereas the
Korteweg-de Vries equation is the central equation to which
Korteweg and De Vries frequently revert in the course of their
paper \cite{ref1}. A simple substitution of the wave velocity into
the continuity equation yields immediately the Korteweg-de Vries
equation in its full glory. However, Boussinesq did not do this,
otherwise it may well be that the history of the long stationary
wave had taken a different course. Pego \cite{ref18} pointed out
that the Korteweg-de Vries equation appeared already in a
footnote on page 360 of Boussinesq's 680 pages vast volume {\it
``M\'emoir sur la th\'eorie des eaux courantes'' } \cite{ref10},
that appeared in 1877, well before the publication of the
Korteweg-de Vries paper in 1895. However, this footnote on the
Korteweg-de Vries equation and also Boussinesq's expression for
the wave velocity are only valid when the wave vanishes at
infinity, while this is not necessary in the theory of Korteweg
and De Vries. Therefore, Boussinesq uses another approach for
treating \textit{steady periodic} waves than Kor\-te\-weg and De
Vries, who presented a unified treatment for steady waves, not
only vanishing at infinity but for waves being periodic as well.
It is not only the equation , but also its applicability that is
important. It seems that this is not always sufficiently realized
or even mentioned in the literature. Darrigol spends in his essay
only one page to the Korteweg-de Vries equation under the heading
{\it ``The so-called Korteweg-de Vries equation'' } \cite{ref17}.
It is only a whim of Tyche, the daughter of Zeus and the
personification of fate, that Zabusky and Kruskal attributed the
names of Korteweg and De Vries to our equation and not that of
Boussinesq, who merits of course the token of priority.

In the following account we present a review of the work by
Boussinesq and Korteweg and De Vries; as to Boussinesq, most of
our attention is directed to his long article in the Journal des
Mathematiques Pures et Appliqu\'ees \cite{ref9}, which is more
accessible than his vast memoir \cite{ref10}. The author, not a
historian, is well aware that he may have overlooked or deleted
important facts, but nevertheless he hopes that this study may
disclose some generally unknown aspects of the early history of
the Korteweg-de Vries equation.

\section{Scott Russell's experiments}

The story of the discovery of the ``Wave of Translation'' by John
Scott Russell in 1834, has been recorded in many books concerning
 ``Soliton'' theory, the more so because Scott Russell's account is
fascinating and even full of emotion, hardly expected in a
scientific paper. Therefore, our account will be rather short and
the interested reader is referred to his ``Report on Waves''
\cite{ref6, ref4, ref12, ref14}. It was in the year of 1834 that
the Scottish naval architect followed on horseback a towboat,
pulled by a pair of horses along the Union Canal, connecting
Edinburgh and Glasgow. However, the boat was suddenly stopped in
its speed - presumably by some obstacle - but not the mass of
water, which it had put in motion. Our engineer perceived a very
peculiar phenomenon: a nice round and smooth wave - a well defined
heap of water - loosened itself from the stern and moved off in
forward direction without changing its form with a speed of about
eight miles an hour and about thirty feet long and one or two feet
in height. He followed the wave on his horse and after a chase of
one or two miles he lost the heap of water in the windings of the
channel \cite{ref6}. Many a physicist would not be inclined to
analyze this phenomenon and leave it as it is, not so Scott
Russell discovering something very peculiar in a seemingly
ordinary event. He designed experiments generating long waves in
long shallow basins filled with a layer of water and he
investigated the phenomenon he had observed. He studied the form
of the waves, their speed of propagation and stability, clearly
perceptible in progressing positive waves, but not in progressing
negative waves. A schematic view of these experiments is shown in
figure 1, which is adapted from Remoissenet \cite{ref21}.

\begin{center}
\includegraphics[width=4in]{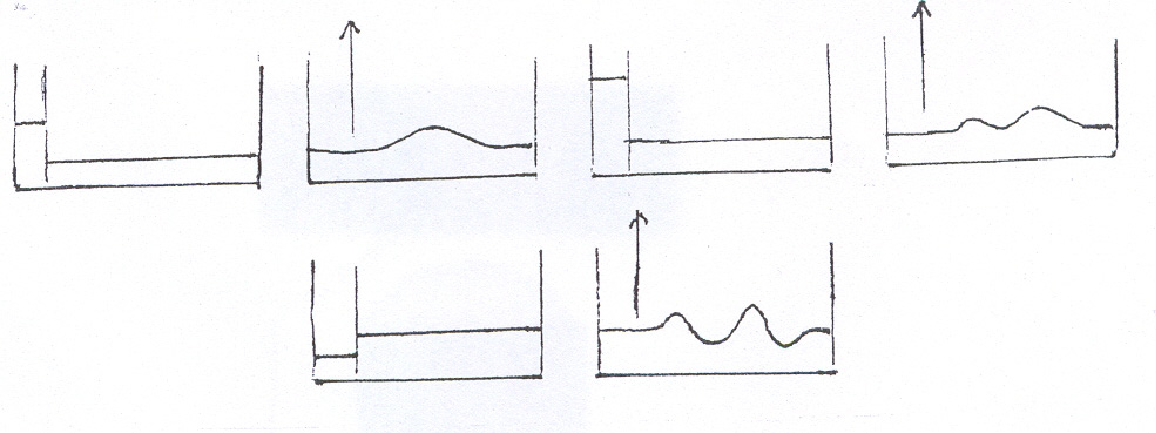}\\
Figure 1: Scott Russell's experiments
\end{center}

For an extensive historical study the reader may consult the
papers by Bullough \cite{ref14} and Darrigol \cite{ref17}. As
mentioned in the introduction, there existed in England and France
a rich tradition in the mathematical description of hydrodynamic
phenomena such as wave motions in fluids. Scott Russell challenged
the mathematical community to prove theoretically the existence of
his solitary wave and ``to give an a priori demonstration a
posteriori'' i.e. to show the possible existence of a stable
solitary wave propagating without change of form.

It is not unusual that new discoveries or new ideas encounter
resistance from established convictions. We take from Rayleigh's
paper ``On Waves'' in the Philosophical Magazine [11, pp 257-279,
1876] the following quotations. Airy, an authority on the subject,
writes in his treatise on ``Tides and Waves'' \cite{ref22}: {\it
``We are not disposed to recognize this wave (discovered by Scott
Russell) as deserving the epithets ``great'' or ``primary'',à and
we conceive that ever since it was known that the theory of
shallow waves of great length was contained in the equation $\frac
{\partial^2 X}{\partial  t^2} = g \kappa \frac {\partial^2 X}{\partial
x^2}$ the theory of the solitary wave has been perfectly well
known''}.
 Further {\it ``Some experiments were made by Mr. Russell on what he calls a negative wave.
 But (we know not why) he appears not to have been satisfied with
 these experiments and had omitted them in his abstract.
 All of the theorems of our IVth section, without exception,
 apply to these as well as to positive waves,
 the sign of the coefficient only being changed''}.
 Probably it was also Airy who expressed for the first time as
 his opinion that long waves in a canal with rectangular cross section
 must necessarily change their form as they advance,
becoming steeper in front and less steep behind and in this he was
supported by Lamb and Busset \cite{ref1, ref22, ref23}. Stokes
believed that the only permanent wave should be basically
sinusoidal, but later on he admitted that he had made a mistake,
see also our section 7.

On the other hand he writes \cite{ref24}: {\it ``It is the opinion
of Mr. Russell that the solitary wave is a phenomenon ``sui
generis'', in no wise deriving its character from the circumstance
of the generation of the wave. His experiments seem to render this
conclusion probable. Should it be correct, the analytical
character of the solitary wave remains to be discussed''}.

The {\it ``a priori demonstration a posteriori''} asked for by
Scott Russell was finally given, first by J. Boussinesq in 1871
[7-10], some time later in 1876, by Lord Rayleigh \cite{ref11} and
in order to remove still existing doubts over the existence of the
solitary wave by G. de Vries \cite{ref25} and by D.J. Korteweg and
G. de Vries in 1895 \cite{ref1}.

In the next section we present first a concise account of the
contribution by Rayleigh, since it is short and leads directly to
the heart of the matter. Moreover, this paper has been of great
influence on the research of Korteweg and De Vries. Consecutively,
we discuss in the other sections the investigations of Boussinesq
and Korteweg-de Vries and we finish with some concluding remarks.

\section{Rayleigh's Solution}

Be given an incompressible irrotational flow in a canal with a
constant rectangular cross-section, fig. 2.

\begin{center}
\includegraphics[width=4in]{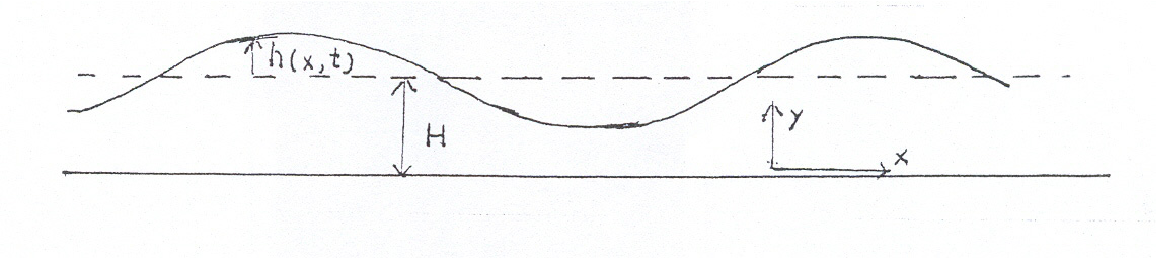}\\
Figure 2: Wavesurface
\end{center}

The coordinates of a fluid particle are given by the coordinates
$x$ and $y$, the undisturbed depth of the canal by the constant
$H$ and the wave surface by $H + h (x,t)$. Another essential
assumption is that the wave length is large in comparison with
$H$.

As has been mentioned in the preceding section, Airy had already
remarked that the theory of shallow waves of great length is
contained in the equation
$$
\frac {\partial^2 X}{\partial t^2}= g\kappa \frac {\partial^2
X}{\partial x^2}
$$
(where $\kappa = H$ and $g$ the constant of gravity). The wave
velocity is $\sqrt{g\kappa}$ , a result generally known since
Lagrange in 1786. However, this result is only valid as a first
order approximation, where $h/H$ may be neglected.

Rayleigh remarked that for this value of the wave velocity the
so-called free surface condition (equilibrium of pressure) is only
satisfied whenever the ratio $h/H$ may be neglected, but if this
is not the case it is impossible to have a wave in still water
with velocity $\sqrt{g\kappa}$ and at the same time propagating
without change of form. In order to cure this discrepancy with
Scott Russell's experimental results, he proposes to look for a
more accurate approximation of the wave velocity (\cite{ref11}, pp
252-253). Rayleigh assumes the existence of a stationary wave,
vanishing at infinity, and by adding to the fluid a yet unknown
constant basic velocity equal and opposite to that of the wave, he
may omit the dependence on time. Since the flow is free of
rotation, and incompressible, there exist a velocity potential
$\phi$ and a stream function $\psi$, both satisfying Laplace's
equation. The horizontal and vertical velocity components are
given by
$$
  u=\frac {\partial \phi}{\partial x}=\frac {\partial\psi}{\partial y}
  \quad\mbox{and}\quad
   v=\frac {\partial \phi}{\partial y}=-\frac {\partial\psi}{\partial x}
$$
and a series expansion gives
\begin{eqnarray}\label{eq3.1}
u&=&\frac {\partial \phi}{\partial x}\,=\, f(x)- \frac{y^2}{2!}
f''(x)+\frac{y^4}{4!} f^{(4)}(x)-\ldots\\
\label{eq3.2} v&=&\frac {\partial \phi}{\partial y}\,=\, -y f'(x)+
\frac{y^3}{3!}f^{(3)}(x)-\ldots
\end{eqnarray}
This expansion is justified because $f(x)$ is, due to the large
wave lengths, a slowly varying function of $x$.

Integration yields the stream function
\begin{equation}\label{eq3.3}
  \psi = y f(x)-\frac{y^3}{3!}f''(x) + \frac{y^5}{5!}f^{(4)}(x)
\end{equation}
constant along stream lines and hence also along the wave surface
$y = H + h(x)$.

Let $p$ be the pressure just below the wave surface, then we have
the relation
\begin{equation}\label{eq3.4}
-2\frac{p-C}{\rho}= 2g(H+h)+ u_s^2+ v_s^2:=\tilde{\omega}
\end{equation}
where $\rho$ is the density, $g$ the constant of gravity and $C$
an integration constant; the suffix $s$ denotes that the values of
$u$ and $v$ are taken at the surface of the wave.

To satisfy in higher approximation the free surface condition, ---
$p$ constant ---, Rayleigh investigates how far it is possible to
make $\tilde{\omega}$ constant by varying $h(x)$ as function of
$x$. Using $u_s^2+ v_s^2 =u_s^2\{1+(\frac {dh}{dx})^2\}$ and
eliminating the unknown function $f$ with the aid of \eqref{eq3.1}
and \eqref{eq3.3} he obtains after a tedious calculation a
differential equation for the wave form $y = H + h(x)$, viz.
\begin{equation}\label{eq3.5}
\tfrac 13 \left(\frac {dy}{dx}\right)^2= 1 + C y+ \frac{
u_0^2+2gH}{u_0^2H^2} y^2- \frac{g}{u_0^2 H^2}y^3
\end{equation}
where $u_0$ is the still unknown constant basic velocity that has
been added to the flow and $C$ is again an integration constant,
(\cite{ref11}, pp 258-259). The cubic expression at the right hand
side vanishes for $y = H$ with $x=\infty$ and for $y = H + h_0$
with $h_0$ the crest of the wave.

Elimination of $C$ yields
\begin{equation}\label{eq3.6}
  u_0= \sqrt{g(H+h_0)}\approx \sqrt{g H}+ \tfrac 12 h_0
  \sqrt{\frac gH}
\end{equation}
which is also the wave velocity, and the equation \eqref{eq3.5}
reduces to
$$
 \left(\frac {dh}{dx}\right)^2+\frac 3{H^3}h^2(h-h_0)=0
$$
with the solution
\begin{equation}\label{eq3.7}
  h(x)=h_0 \,{\mathrm{sech}}^2\left(\sqrt{\frac {3h_0}{4H^3}}x\right).
\end{equation}

This formula represents the ``heap of water'' (the ``Great Wave'')
with the right wave velocity \eqref{eq3.6}, as experimentally
observed by Scott Russell and so he was finally vindicated after
fourty years of much discussion in England.

Rayleigh finishes his article with the remark: {\it ``I have
lately seen a memoir by M. Boussinesq, Comptes Rendus, Vol. LXXII,
in which is contained a theory of the solitary wave very similar
to that of this paper. So as far as our results are common, the
credit of priority belongs of course to M. Boussinesq''}.

\section{The Equations of Boussinesq}

{}From Miles we learn in his essay \cite{ref16} that Boussinesq
(1842-1929) received his doctorate from the Facult\'e des
Sciences, Paris, in 1867, occupied chairs at Lille from 1873 to
1885, and at the Sorbonne from 1885 to 1896. He made significant
contributions to hydrodynamics and the theories of elasticity,
light and heat. He wrote several papers on nonlinear dispersive
waves \cite{ref7, ref8, ref9} and a voluminous ``m\'emoire",
entitled ``Essai sur la th\'eorie des eaux courantes", presented
to the Acad\'emie des Sciences in 1877, Vol. XXIII (\cite{ref10},
pp 1-680).

The most accessible publication is his article in the Journal de
Math\'e\-ma\-ti\-ques Pures et Appliqu\'ees in 1872 \cite{ref9}.
It has the verbose title {\it ``Th\'eorie des ondes et des remous
qui se propagent le long d'un canal rectangulaire horizontal, en
communiquant au liquide continu dans ce canal des vitesses
sensiblement pareilles de la surface au fond''}. This paper
subsumes the short Comptes Rendus \cite{ref7, ref8}, whereas the
monograph \cite{ref10} gives also much information less relevant
for this exposition concerning the ``Great Wave".

He considers, in the same way as Rayleigh, long waves in a shallow
canal with rectangular cross section; the fluid is supposed
incompressible and rotation free, while friction, also along the
boundaries of the canal, is neglected. Distinct from Rayleigh's
article, Boussinesq introduces also a time variable, essential for
the description of a dynamic phenomenon, and the coordinates of a
fluid particle at time $t$ are denoted by $(x,y) = (x(t), y(t))$,
see fig. 2.

Be $p$ the pressure in the fluid, $\rho$ its density and $(u,v)$
the velocity vector. The height of the water in equilibrium is
again denoted by the constant $H$ and the wave surface by the
function $y = H + h(x,t)$. The wave length is supposed to be large
and the amplitude $h$ of the wave small in comparison with $H$ and
vanishing for $x\to \pm \infty$.

Boussinesq's exposition starts along the same lines as in the
theory of Ray\-leigh. The main ingredients are a series
development into powers of $y$, similar as in \eqref{eq3.1} --
\eqref{eq3.3}:
\begin{equation}\label{eq4.1}
 \phi= f- \frac{y^2}{2!}
f''+\frac{y^4}{4!} f^{(4)}-\ldots
\end{equation}
with the as yet unknown function $f = f(x,t)$ and valid for $0 <
y< H + h(x,t)$.

Further, he uses the free surface condition
\begin{equation}\label{eq4.2}
  gh+\tfrac 12 (u_s^2+v_s^2)+ \frac{\partial \phi_s}{\partial t} +
  \chi(t)=0
\end{equation}
where $\chi(t)$ is an arbitrary function and where the suffix $s$
refers to the wave surface. Under the assumption that the
potential $\phi$ and its derivatives with respect to $x$, $y$ and $t$
vanish for $x \to\pm \infty $, the function $\chi(t)$ may be
omitted. We shall see that this assumption is very essential in
Boussinesq's theory and it is used again and again in his paper,
see also section 7.

A second boundary condition follows from the kinematic equation
\begin{equation}\label{eq4.3}
  v_s=\frac{dh}{dt} = \frac{\partial h}{\partial t}+u_s\frac{\partial h}{\partial x}
\end{equation}

Substitution of the series expansions of the potential and the
velocity components into \eqref{eq4.2} and \eqref{eq4.3} results
into two equations, containing $h(x,t)$ and $f(x,t)$. Elimination
of $f(x,t)$ gives in a first approximation, with $h(x,t)$ small in
comparison with $H$, the wave equation of Lagrange:
$$
\frac{\partial^2 h}{\partial t^2}= gH\frac{\partial^2 h}{\partial
x^2}
$$
 A second higher approximation yields
the well known equation of Boussinesq
\begin{equation}\label{eq4.4}
  \frac{\partial^2 h}{\partial t^2}=gH\frac{\partial^2 h}{\partial
x^2}+g H\frac{\partial^2 }{\partial x^2}\left[ \frac{3h^2}{2H}+
\frac{H^2}{3}\frac{\partial^2 h}{\partial x^2}\right]
\end{equation}
This equation may be simplified by restricting the theory to waves
propagating into only one direction, say the positive $x$-axis.
Denoting the wave velocity by $\omega(x,t)$ and using the
conservation of mass
\begin{equation}\label{eq4.5}
\frac{\partial h}{\partial t}+\frac{\partial}{\partial x}(\omega
h) =\frac{\partial h}{\partial t}+\omega\frac{\partial h}{\partial
x}+h\frac{\partial \omega}{\partial x}=0
\end{equation}
we obtain after substitution into \eqref{eq4.4} and integration
with respect to $x$
\begin{equation}\label{eq4.6}
\frac{\partial }{\partial t}(\omega h) + g H \frac{\partial
}{\partial x}(h+\frac 32\frac {h^2}{H} +
\frac{H^2}{3}\frac{\partial^2 h}{\partial x^2})=0
\end{equation}

To make progress it is desirable to have an explicit expression
for $\omega(x,t)$, because substitution into \eqref{eq4.5} or
\eqref{eq4.6} gives a differential equation for the wave surface
$h(x,t)$. To this end Boussinesq introduces without a clear
motivation the function
\begin{equation}\label{eq4.7}
\psi(x,t)=h\cdot(\omega-\sqrt{g H})-
\frac{\sqrt{gH}}{2}\left(\frac 32\frac {h^2}{H} +
\frac{H^2}{3}\frac{\partial^2 h}{\partial x^2}\right).
\end{equation}

Differentiation of this expression with respect to $t$ and
replacing $\frac{\partial }{\partial t}$ by $-\sqrt{g
H}\frac{\partial }{\partial x}$, which does not disturb the order
of approximation in the second term, yields
$$
\frac{\partial \psi}{\partial t}= \frac{\partial }{\partial
t}(\omega h)-\sqrt{g H}\frac{\partial h}{\partial t} +\frac {g
H}{2}\frac{\partial }{\partial x}\left(\frac 32\frac {h^2}{H} +
\frac{H^2}{3}\frac{\partial^2 h}{\partial x^2}\right)
$$

Substitution of \eqref{eq4.6} and \eqref{eq4.5} gives
$\frac{\partial \psi}{\partial t}=\sqrt{g H}\frac{\partial
\psi}{\partial x}$ , from which it follows that $\psi\equiv 0$ for
a wave propagating in the positive $x$-direction. Hence,
Boussinesq obtains from \eqref{eq4.7} the important result
\begin{equation}\label{eq4.8}
\omega(x,t) =\sqrt{g H} + \sqrt{g H}\left(\frac {3h}{4H} +
\frac{H^2}{6h}\frac{\partial^2 h}{\partial x^2}\right)
\end{equation}

{}From \eqref{eq4.5} and \eqref{eq4.8} one may obtain the
differential equations
$$
h\frac{d h}{d t}= h \left( \frac{\partial h}{\partial
t}+\omega\frac{\partial h}{\partial x}\right)= - h^2
\frac{\partial \omega}{\partial x}= -\frac{\partial }{\partial
x}(h^2\omega)+ 2 \omega h \frac{\partial h}{\partial x}
$$
 and
\begin{equation}\label{eq4.9}
\frac{d h}{d t}= -\frac 14 \sqrt{\frac gH}\frac 1h \frac{\partial
}{\partial x}\left[ h^3\left\{ 1+\frac 23 H^3 \frac 1h
\left(\frac{\partial }{\partial x}\frac 1h \frac{\partial
h}{\partial x}\right)\right\}\right]
\end{equation}
or after passing to the new variable $h\, dx = -d\sigma$
\begin{equation}\label{eq4.10}
\frac{d h}{d t}=\frac 14 \sqrt{\frac gH}\frac{\partial }{\partial
\sigma}\left[ h^3\left\{ 1+ \frac 23 H^3 \frac{\partial^2 h
}{\partial \sigma^2}\right\}\right].
\end{equation}

It appears from \eqref{eq4.8} that the wave velocity differs from
point to point at the wave surface and so it is expected that the
wave should change its form during its course, which is one of the
main issues in pursuing Scott Russell's experiments. However, the
wave will only be stationary whenever $\omega(x,t)$ is constant.
It is appropriate to make here some comments:

1.  The introduction of the function $\psi$ in \eqref{eq4.7} is
not a priori clear and motivated. With the aid of \eqref{eq4.5} we
may write instead of \eqref{eq4.6}
$$
\frac{\partial }{\partial t}\left\{ h(\omega-\sqrt{gH})\right\}
-\sqrt{gH}\frac{\partial }{\partial x}\left\{
h(\omega-\sqrt{gH})\right\}+ gH\frac{\partial }{\partial x}\left(
\frac{3h^2}{ 2H} + \frac{H^2}{3} \frac{\partial^2 h }{\partial
x^2}\right)=0.
$$
Replacing $\frac{\partial }{\partial t}$ by $-\sqrt{gH}
\frac{\partial }{\partial x}$ and integrating with respect to $x$,
we get the result \eqref{eq4.8}; remember $h(x,t)$ and its
derivatives vanish for  $x \to \pm \infty$.

2. Another proof of \eqref{eq4.8} was given many years later in
1885 by St.~Venant \cite{ref19}. It may be that he was not
satisfied with Boussinesq's derivation. He applied another
approach, using the mean value of the horizontal component of the
velocity vector
$$U(x,t) = \frac{1}{H+h(x,t)} \int_0^{H+h(x,t)}u\, dy
$$

3.  The equations \eqref{eq4.9} and \eqref{eq4.10} contain
implicitly the wave velocity $\omega$. This dependence on $\omega$
can simply be eliminated by the substitution of \eqref{eq4.8} into
\eqref{eq4.5}. If Boussinesq had carried out this small operation
he had obtained the Korteweg-de Vries equation long ``avant la
lettre", viz.
\begin{equation}\label{eq4.11}
\frac{\partial h}{\partial t}+ \sqrt{\frac gH }\frac
32\frac{\partial }{\partial x}\left(\frac 23H h + \frac 12 h^2 +
\frac{H^3}9 \frac{\partial^2 h }{\partial x^2}\right)=0
\end{equation}
valid for waves vanishing at infinity. This equation does not
differ essentially from the Korteweg-de Vries equation as
presented in the Korteweg-de Vries paper in the Philosophical
Magazine; the coordinates in \eqref{eq4.11} refer to a fixed
$(x,t)$ frame, whereas Korteweg and De Vries used a moving frame,
see also next section.

4.  As mentioned in the Introduction of this paper, R.~Pego
\cite{ref18} and also O.~Darrigol (\cite{ref17}, pp 47) have
discovered in a footnote on page 360 of the ``Essai sur la
th\'eorie des eaux courantes'', that Boussinesq had found already
in 1876 the Korteweg-de Vries equation, mind without recourse to
the expression \eqref{eq4.8} for the wave velocity \cite{ref10}.
 In fact he used instead of \eqref{eq4.7} the function
$$
\psi_1(x,t) = \frac{\partial h}{\partial t}+ \sqrt{gH }
\frac{\partial h}{\partial x}+ \frac12 \sqrt{gH} \frac{\partial
}{\partial x}\left(\frac {3h^2}{2 H} + \frac {H^2}{3}
\frac{\partial^2 h }{\partial x^2}\right)= -\frac{\partial
\psi}{\partial x}
$$
and using $\psi \equiv 0$ he gets equation \eqref{eq4.11}.

Consecutively, the wave velocity $\omega(x,t)$ may be determined
by the integration of \eqref{eq4.5}, i.e.
$$\omega(x,t) = \frac{1}{h(x,t)} \int_{-\infty}^{x}\left(-\frac{\partial h}{\partial
t}\right) dx
$$

5.  The bidirectional Boussinesq equation \eqref{eq4.4} can be
factorized as
\begin{equation}\label{eq4.12}
  \left(\frac{\partial }{\partial t}- \sqrt{gH} \frac{\partial }{\partial
  x}\right)
  \left\{\frac{\partial h}{\partial t}+\sqrt{\frac gH}\frac 32 \frac{\partial }{\partial x}
  \left(\frac 23 H h + \frac 12 h^2 + \frac{ H^3}{9} \frac{\partial^2h}{\partial
  x^2}\right)\right\}=0
\end{equation}
from which it immediately follows that the unidirectional
Korteweg-de Vries equation \eqref{eq4.11} is contained in
Boussinesq's equation \eqref{eq4.4}.

\section{The Korteweg-de Vries Equation}

We start with a few biographical data of Korteweg and De Vries
\cite{ref20}. Diederik Johannes Korteweg (1848-1941) received in
1878 his doctorate at the University of Amsterdam, after defending
his thesis on the propagation of waves in elastic tubes. His
supervisor was J.D. van der Waals, renowned for his equation of
state and the continuity of the gas and fluid phases. Korteweg
occupied the chair of mathematics, mechanics and astronomy at the
University of Amsterdam from 1881 to 1918; he published several
papers on mathematics, classical mechanics, fluid mechanics and
thermodynamics. We mention in particular his investigations on the
properties of ``folded'' surfaces in the neighbourhood of singular
points, work related to that of Van der Waals \cite{ref26}.
Another scientific achievement is the edition of the ``Oevres
Compl\`etes'' of Christiaan Huygens and Korteweg was the principal
leader of this project in the period 1911-1927.

He inspired many young mathematicians who wrote their thesis under
his supervision, among others Gustav de Vries and the famous
L.E.J. Brouwer. Korteweg had a great influence on academic life in
the Netherlands as appears from his leadership in several academic
institutions. The thesis of Gustav de Vries, entitled ``Bijdrage
tot de Kennis der Lange Golven'' \cite{ref25} was published in
1894 and the paper in the Philosophical Magazine of 1895 is an
excerption of this thesis. De Vries published papers on cyclones
in 1896 and 1897 and two papers ``Calculus Rationum'' in the
proceedings of the Royal Dutch Academy of Arts and Sciences in
1912. He teached mathematics at a gymnasium in Alkmaar and at a
secondary school in Haarlem.

The author received from the grandsons of De Vries a copy of the
scientific correspondence between Korteweg and De Vries. From this
we know that Korteweg advised De Vries to study Rayleigh's method
of the series expansion which has been explained in section 3 of
this paper. He also suggested to include capillarity and to
investigate long periodic waves.

Korteweg and De Vries start their article with the time dependent
Rayleigh expansions
\begin{eqnarray}\label{eq5.1}
  u(x,t)&=& f(x,t)- \frac{y^2}{2!}
f''(x,t)+\frac{y^4}{4!} f^{(4)}(x,t)-\ldots\\
\label{eq5.2} v(x,t)&=& -y f'(x,t)+
\frac{y^3}{3!}f^{(3)}(x,t)-\ldots
\end{eqnarray}

The effect of the surface tension in the free surface condition
amounts to an extra term in \eqref{eq4.2}:
\begin{equation}\label{eq5.3}
  gh+\tfrac 12 (u_s^2+v_s^2)+ \frac{\partial \phi_s}{\partial t} +
  \chi(t)=\frac T\rho\frac{\partial^2 y_s}{\partial x^2}
\end{equation}
where $T$ is the surface tension and $\chi(t)$ again the arbitrary
function depending only on time.

This arbitrary function is eliminated by Boussinesq with the aid
of the assumption that $f$ and its derivations vanish for $x \to
\pm\infty$. Korteweg and De Vries drop this crucial assumption and
$\chi(t)$ is eliminated by simply dif\-fe\-ren\-tia\-ting \eqref{eq5.3}
with respect to $x$. It is now already noted that {\it periodic}
waves are not a priori excluded from further discussion, this in
contrast to Boussinesq, who used a different approach in his
discussion of periodic waves, see section 7.1. Differentiating
\eqref{eq5.3} with respect to $x$, the free surface condition
\eqref{eq5.3} becomes
\begin{equation}\label{eq5.4}
  g\frac{\partial h}{\partial x}+ u_s\frac{\partial u_s}{\partial
  x}+ v_s\frac{\partial v_s}{\partial x}+ \frac{\partial^2 \phi_s}{\partial
  t\partial x}- \frac T\rho \frac{\partial^3 y_s}{\partial x^3}=0
\end{equation}
Besides this we need again the kinematic condition \eqref{eq4.3}
\begin{equation}\label{eq5.5}
v_s=\frac{\partial h}{\partial t}+ u_s\frac{\partial h}{\partial
x}.
\end{equation}

Korteweg and De Vries put $y_s = H + h(x,t)$ and $f(x,t) = q_0 +
\beta(x,t)$ with $q_0$ an as yet undetermined constant velocity.
Substitution of \eqref{eq5.1} and \eqref{eq5.2} into \eqref{eq5.4}
and \eqref{eq5.5} gives in a first order approximation for $h$
small and for a wave, progressing in the positive $x$-direction,
the expression
$$h = h(x - (q_0 + \sqrt{gH})t).$$

Adding to the flow a velocity $q_0 = - \sqrt{gH}$, we obtain the
Lagrange steady wave solution with 
$$\frac{\partial h}{\partial
t}=0,\qquad \frac{\partial \beta}{\partial t}=0,$$ 
and 
$$\frac{\partial
\beta}{\partial x}= -\frac {q_0}{H} \frac{\partial h}{\partial x}=
- \frac{g}{q_0}\frac{\partial h}{\partial x}\mbox{ or }\beta =
-\frac{g}{q_0} (h + a),$$ 
where $a$ is an undetermined constant.

The next approximation is obtained by $$ f(x,t) = q_0
-\frac{g}{q_0}(h(x,t) + \alpha + \gamma(x,t))
$$
 with $\gamma$ small in
comparison with $h$ and $a$. Substitution into \eqref{eq5.4} and
\eqref{eq5.5} yields two equations for $h$ and $\gamma$ and
elimination of $\gamma(x,t)$ gives finally the Korteweg-de Vries
equation as it appeared for the first time in the thesis of De
Vries \cite{ref25}:
\begin{equation}\label{eq5.6}
\frac{\partial h}{\partial t}= \frac 32 \frac {g}{q_0}
\frac{\partial }{\partial x}\left( \frac 12 h^2 + \frac 23
\alpha h + \frac 13 \sigma \frac{\partial^2 h}{\partial
x^2}\right)
\end{equation}
with $\sigma = \frac 13 H^3 - \frac{TH}{\rho g}$.

The addition of the velocity $q_0 - \frac{g}{q_0}\,\alpha =
-\sqrt{gH} +\sqrt{g/H} \alpha$ to the flow may also be obtained by
a transformation of the fixed $(x,y)$ coordinate system to the
moving frame
\begin{equation}\label{eq5.7}
\xi = x- (\sqrt{gH} -\sqrt{\frac gH} \alpha)t, \quad \tau=t.
\end{equation}

Hence, in this moving frame and forgetting about the added
velocity, we get
\begin{equation}\label{eq5.8}
\frac{\partial h}{\partial \tau}+\frac 32 \sqrt{\frac gH}
\frac{\partial }{\partial \xi}\left( \frac 12 h^2 + \frac 23
\alpha h + \frac 13 \sigma \frac{\partial^2 h}{\partial
\xi^2}\right)=0.
\end{equation}
This equation with $T = 0$ is equivalent with Boussinesq's result
\eqref{eq4.11} and may be obtained by substititution of
\eqref{eq5.7} into \eqref{eq4.11}.

\section{The Solitary Great Wave}

Because of the equivalence of the differential equations
\eqref{eq4.11} and \eqref{eq5.8} for the surface of a wave with
amplitude, vanishing at infinity, it is evident that the theory of
Boussinesq leads to the same results as that of Korteweg and De
Vries, if capillarity is neglected. A necessary and sufficient
condition for the existence of a steady wave is a constant uniform
wave velocity in all points of the wave surface.

\subsection{The Solitary Steady Wave in the Theory of Boussinesq}

It follows from \eqref{eq4.8} that the wave is stationary if
$\omega = \sqrt{gH} + \frac 12 \sqrt{g/H} h_1$, where $h_1$ is
some as yet unknown constant, independent of $x$ and $t$.
Therefore,
\begin{equation}\label{eq6.1}
\frac{\partial^2 h}{\partial x^2}= \frac{3h }{ 2H^3}(2h_1-3h)
\end{equation}
Integration with $h \to 0$  and $\frac{\partial h}{\partial x}\to
0$ for $x \to\pm\infty$ gives
\begin{equation}\label{eq6.2}
h(x,t) = h_1
{\mathrm{sech}}^2\left\{\sqrt{\frac{3h_1}{4H^3}}(x-\omega
t)\right\}.
\end{equation}

It follows $h_1 \geq h(x,t)$  and $h_1$ is the crest of the wave.
The wave velocity
\begin{equation}\label{eq6.3}
  \omega= \sqrt{gH} + \frac 12 \sqrt{\frac gH}h_1
\end{equation}
contains a correction of Lagrange's result with $\omega =
\sqrt{gH}$ and was already experimentally verified in 1844 by
Scott Russell \cite{ref6}. The expressions \eqref{eq6.2} and
\eqref{eq6.3} were also obtained by Rayleigh, see \eqref{eq3.6}
and \eqref{eq3.7}. It appears that the wave velocity is the larger
the higher the crest of the wave and this means that in the case
of several separate solitary waves of the form \eqref{eq6.2} the
higher waves will overtake the lower ones whenever the higher
waves were initially behind the lower waves. It is known that this
occurs without change of form, however, there is only a phase
shift.

The solitary waves behave like a row of rolling marbles, where the
faster marbles carry over their impuls to the slower marbles. They
were coined by Zabusky and Kruskal \cite{ref3} as ``solitons'' to
indicate their particle like properties. For an explicit
calculation of this behaviour the interested reader may consult
ref. \cite{ref13}, part II, 3.5.

\subsection{The Solitary Steady Wave in the Theory of Korteweg
and De Vries}

Korteweg and De Vries do not have at their disposal an explicit
expression for the wave velocity. However, for a steady wave in
the moving $(\xi,\tau)$ frame \eqref{eq5.7} one has
$\frac{\partial h}{\partial \tau} = 0$ and so by \eqref{eq5.8}
\begin{equation}\label{eq6.4}
\frac{d}{d\xi}\left( \frac 12 h^2 +\frac 23 \alpha h + \frac 13\sigma
\frac{d^2h}{d\xi^2}\right)=0
\end{equation}
where $\alpha$ is the still unknown correction of the wave
velocity. Integration under the assumption $h$, $\frac
{dh}{d\xi}$, $\frac {d^2h}{d\xi^2}\to  0$ for $\xi\to \pm\infty$
results into
$$\frac{dh}{ d\xi} = \pm \sqrt{\frac{-h^2(h+2\alpha)}{\sigma}}.$$
There are two distinct cases, $\sigma > 0$ and $\sigma < 0$; we
restrict our calculation to the case $\sigma > 0$ since the other
case can be treated similarly. If $\sigma > 0$ then the constant
$2 \alpha$ is negative and taking $2 \alpha = - h_2$ one gets
\begin{equation}\label{eq6.5}
  h(\xi)=
  h_2{\mathrm{sech}}^2\left(\sqrt{\frac{h_2}{4\sigma}}\xi\right),
\end{equation}
with  $h_2 > 0$ the crest of the wave. When the surface tension
$T$ is neglected, the parameter $\sigma$ reduces to $\sigma =
\frac 13 H^3$ and \eqref{eq6.5} is in agreement with Boussinesq's
result \eqref{eq6.2}. Also the wave velocity is conform
\eqref{eq6.3} because we have by \eqref{eq5.7}
\begin{equation}\label{eq6.6}
  \omega = \sqrt{gH} - \sqrt{\frac gH}\alpha = \sqrt{gH}+\frac
  12\sqrt{\frac gH} h_2
\end{equation}

Korteweg and De Vries consider also solitary steady waves with a
negative amplitude, possible for $\sigma < 0$, i.e. $H<\sqrt{\frac
{3T}{\rho g}}$ ; this quantity equals approximately $\frac 12$ cm
for water.

\section{Periodic Stationary Waves}

\subsection{ Boussinesq's Theory}

Boussinesq investigates in his M\'emoire {\it ``Essai sur la
theorie des eaux courrantes''} la {\it Forme la plus g\'en\'erale
des intumescences propag\'ees le long d'un canal horizontal et
rectangulaire, qui avancent sans se d\'eformer''} (\cite{ref10},
pp 390-396). By now, he does not have a formula for the wave
velocity as in \eqref{eq4.8}, because the condition $h(x,t) \to 0$
for $x\to \pm\infty$ does no longer hold in the case of periodic
waves and so an easy evaluation of the wave form by setting $\omega$
constant is no longer possible. He introduces the mean horizontal
velocity
\begin{equation}\label{eq7.1}
  U(x,t) = \frac{1}{H+h(x,t)} \int_0^{H+h(x,t)}u\, dy
\end{equation}
which satisfies the relation
$$\frac{\partial U}{\partial t}+ \frac{\partial }{\partial
x}\left( gh + \frac 12 U^2 + \frac H3 \frac{\partial^2
h}{\partial t^2}\right)=0
$$
see \cite{ref10}, 276, pp 390-391,  or the paper by St.~Venant
\cite{ref19} where a rather short derivation is presented.

For a steady wave we have $\frac{\partial h}{\partial t}= -
\omega\frac{\partial h}{\partial x}$, $\frac{\partial^2
h}{\partial t^2}= \omega^2\frac{\partial^2 h}{\partial x^2}$ and
$\frac{\partial U}{\partial t}= -\omega \frac{\partial U}{\partial
x}$, so the latter equation becomes
\begin{equation}\label{eq7.2}
  -\omega U + g (h-H ) + \frac {U^2}{2}+ \frac {H\omega^2}{3}\frac{\partial^2 h}{\partial
  x^2}= constant:=\frac{\omega^2}{2H^2} c',
\end{equation}
with $\omega$ the unknown constant wave velocity and
$\frac{\omega^2}{2H^2} c'$ the constant of integration.

The mean velocity $U$ is eliminated with the aid of the
conservation of flux in a reference system bound to the wave, i.e.
$$ \int_0^{H+h(x,t)} u dy = \omega h(x,t)$$
or
\begin{equation}\label{eq7.3}
  U=\frac{\omega h } {H+h}\approx \frac {\omega h}{H}\left(1-\frac hH\right).
\end{equation}

Substitution into \eqref{eq7.2} and omitting terms of order $h^3$
and higher we get
$$2\frac{\partial^2 h}{\partial x^2}= -\frac 3{H^3} \{ 3 h^2 - 2 H
h \left(1 - \frac {g H }{\omega^2}\right)- c'\},$$ and by
multiplication with $\frac{\partial h}{\partial x}$ and
integration with respect to $x$
\begin{equation}\label{eq7.4}
\left(\frac{\partial h}{\partial x}\right)^2= -\frac{3}{H^3}\{ h^3
- H h^2\left(1 - \frac {g H }{\omega^2}\right)- c' h - c''\},
\end{equation}
with $c''$  the integration constant.

\begin{center}
\includegraphics[width=3in]{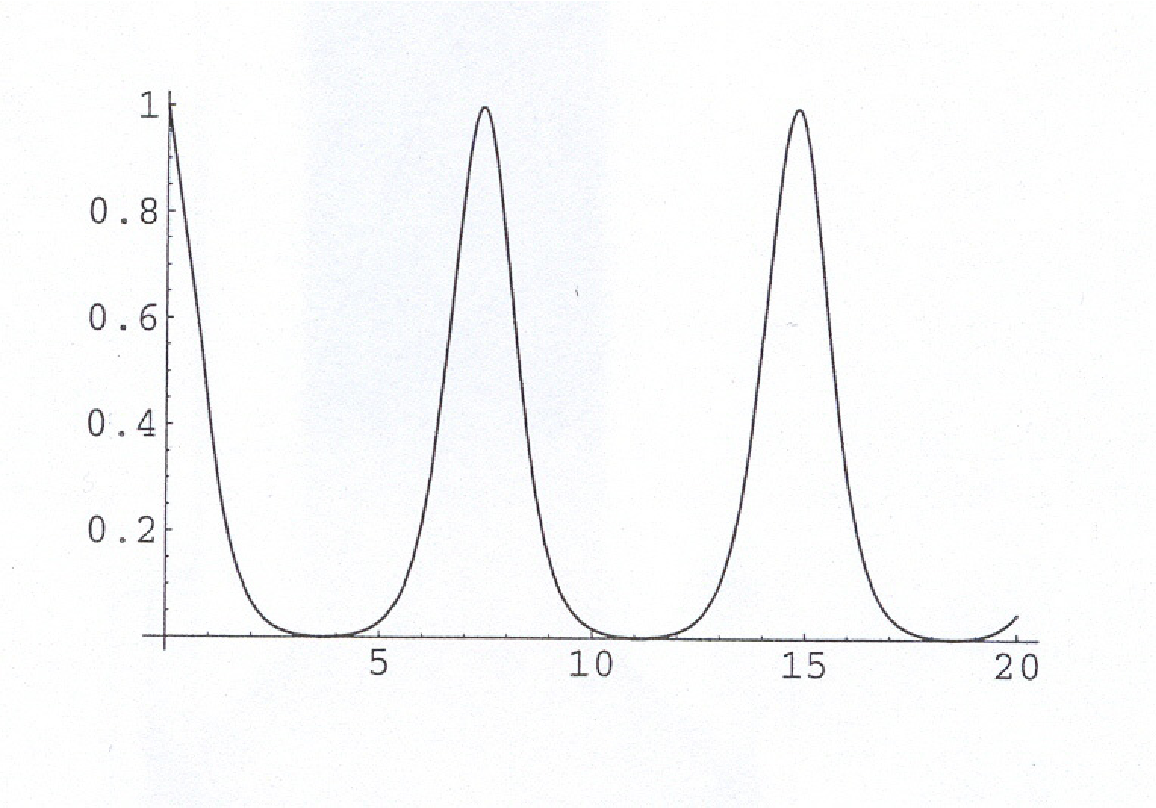}\\
Figure 3: Cnoidal Wave
\end{center}

Suppose now that the minimum value $h_0$ of $h$, which is
necessarily negative, is assumed at $x = 0$ and the wave surface
represented by $\tilde{h}(x,t) =h(x,t)-h_0$, then $\tilde{h}\geq
0$ and $\tilde{h}=0$ at x = 0. This shift results into
\begin{equation}\label{eq7.5}
\left(\frac{\partial \tilde{h}}{\partial x}\right)^2=
-\frac{3}{H^3}\{ \tilde{h}^3 +a_2\tilde{h}^2+ a_1\tilde{h} +
a_0\},
\end{equation}
where $a_0$, $a_1$ and $a_2$ are certain constants, which due to
the condition $\tilde{h}= 0$ and $\frac{\partial
\tilde{h}}{\partial x} = 0$ at $x = 0$ satisfy $a_0 = 0$ and $a_1
< 0$. Therefore, \eqref{eq7.5} may be written as
\begin{equation}\label{eq7.6}
\left(\frac{\partial \tilde{h}}{\partial x}\right)^2=
\frac{3}{H^3}\{ (\tilde{h}+k)\tilde{h}(l-\tilde{h})\},
\end{equation}
with $k$ and $l$ positive constants.

Integration leads to a Jacobian elliptic function, but Boussinesq
recommends to use Newton's binomium for a series expansion of the
variable $x$ into powers of $\tilde{h}$.

It follows from \eqref{eq7.6} that the solitary steady wave is
only a particular case, resulting from \eqref{eq7.6} for $k \to
0$. Finally, one may substitute $h = \tilde{h} + h_0$ in
\eqref{eq7.4}, compare the coefficients of $\tilde{h}^2$ in
\eqref{eq7.4} and \eqref{eq7.6} and one obtains after a simple
evaluation the following expression for the wave velocity
\begin{equation}\label{eq7.7}
\omega^2= g \{ H+(l-k)\}.
\end{equation}

\subsection{The Theory of Korteweg and De Vries}

Since the Korteweg-de Vries equation  \eqref{eq5.8}, as derived
in section 5, may also be applied to waves, not necessarily
vanishing for $x \to \pm\infty$, the equation for the amplitude $h
(\xi)$ of a steady wave is given by
\begin{equation}\label{eq7.8}
\frac{d}{d\xi} \left( \frac 12 h^2 + \frac 23 \alpha h+ \frac
13 \sigma \frac{d^2 h}{d\xi^2}\right)=0.
\end{equation}
Integrating this expression two times one obtains
$$c_1+ \frac 12 h^2+ \frac 23\alpha h + \frac 13 \sigma\frac{d^2 h}{d\xi^2}=0
$$
and
\begin{equation}\label{eq7.9}
c_2+6 c_1 h + h^3 + 2 \alpha h^2 +\sigma \left(\frac{d
h}{d\xi}\right)^2=0,
\end{equation}
with $c_1$ and $c_2$ the constants of integration.

The wave surface may be defined as $y = H_0 + \tilde{h}(\xi)$ with
$H_0$ the minimum value of $y$, $\tilde{h}(\xi) \geq 0$ and
$\tilde{h}(0) = 0$.

It follows that $\frac{d\tilde{h}}{d\xi} = 0$ and
$\frac{d^2\tilde{h}}{d\xi^2}>0$ for $\tilde{h} = 0$ and so $c_2 =
0$ and $c_1 < 0$ under the assumption $\sigma > 0$. Consequently,
the equation $\mu^2 + 2 \alpha \mu + 6 c_1 = 0$ has a positive
root $l$ and a negative root $-k$ and \eqref{eq7.9} reads
\begin{equation}\label{eq7.10}
\frac{d\tilde{h}}{d\xi}= \pm \sqrt{\frac
1\sigma(\tilde{h}+k)\tilde{h}(l-\tilde{h})},
\end{equation}
which for $T = 0$ is the same as \eqref{eq7.6}.

With the aid of the substitution $\tilde{h}= l\cos^2\chi$ Korteweg
and De Vries obtain the well-known periodical ``cnoidal'' wave
\begin{equation}\label{eq7.11}
 \tilde{h}(\xi)= l\,{\mathrm{cn}}^2\left( \sqrt{\frac{l+k}{4\sigma}}
  \xi\right),
\end{equation}
where ${\mathrm{cn}}$ denotes one of the Jacobian elliptic
functions with modulus $M = \sqrt{\frac{l}{l+k}}$, period $$4K = 4
\int_0^1(1-t^2)^{-\frac 12} (1-M^2 t^2)^{-\frac 12}dt$$ and wave
length $4K\sqrt{\frac{\sigma}{l+k}}$.

This wave length becomes infinitely large for $k \to 0$ and the
result is the solitary steady wave \eqref{eq6.5}. However, one
gets for large values of $k$, i.e. for small values of $M$, the
sinusoidal wave
$$\tilde{h}(\xi)= l \cos^2\chi = l \cos^2\left( \sqrt{
\frac{l+k}{4\sigma}}\xi\right)$$ with decreasing wave length for
increasing $k$. This agrees with a result of Stokes \cite{ref24}
and in this case $\tilde{h}(\xi)$ may be expanded in a Fourier series;
this may be the reason why Stokes at first believed that the only
permanent wave should be of sinusoidal type.

The approach of Korteweg and De Vries as given here is in
particular attractive, because of the central role of their
equation \eqref{eq5.8} to which they frequently revert in the
development of their theory.

\section{The Stability of the Stationary Solitary Wave}

It follows from the $(x,t)$ dependence of the wave velocity
$\omega$, \eqref{eq4.8}, that wave propagation involves in general
a change of form, but by definition this does not occur in the
case of a steady wave and so the question arises why the steady
solitary wave is stable and an exception to the rule. For the
possible existence of the steady wave a further investigation is
required, in particular with regard to the ``parameters''
determining the stable behaviour. This has been carried out by
Boussinesq and Korteweg-de Vries in quite different ways. The
presence of the non-linear term $h \frac{\partial h}{\partial x}$
and the dispersion term $\frac 19 H^3\frac{\partial h^3}{\partial
x^3}$ is already an indication for a possible balance, furthering
the stability of the wave.

\subsection{ Stability in the Theory of Korteweg-de Vries }

The authors consider a wave form close to that of the steady
solitary wave
\begin{equation}\label{eq8.1}
h(\xi) = \bar{h} \,{\mathrm{sech}}^2(p \xi)
\end{equation}
where $\bar{h}$ and $p$ are as yet arbitrary constants with $p$
near $\sqrt{\frac {\bar{h}}{4 \sigma}}$ (see \eqref{eq6.5}.

The deformation of this wave is determined by the equation
\eqref{eq5.8} and substitution of \eqref{eq8.1} gives an equation
for the evolution of the surface of the wave, given by
$y=h(\xi,\tau)$:
\begin{equation}\label{eq8.2}
\frac{\partial h}{\partial \tau} = 3 \sqrt{\frac gH}\bar{h} p
(4\sigma p^2- \bar{h}) \left\{- {\mathrm{sech}}^2(p \xi) + \frac
23 \frac{\alpha + 2 \sigma p^2}{4\sigma p^2-\bar{h}}\right\}
{\mathrm{sech}}^2(p\xi) \tanh(p\xi).
\end{equation}
Taking $\alpha =4 \sigma p^2- \frac 32 \bar{h}$ this equation
becomes
\begin{equation}\label{eq8.3}
\frac{\partial h}{\partial \tau} = 3 \sqrt{\frac gH}\bar{h} p
(4\sigma p^2- \bar{h}) {\mathrm{sech}}^2(p\xi) \tanh^3(p\xi).
\end{equation}

The choice $p = \sqrt{\frac {\bar{h}}{4\sigma}}$  and thus $\alpha
= -\frac 12 \bar{h}$ results into $\frac{\partial h}{\partial
\tau}=0$ and we get the steady wave \eqref{eq6.5}.

A numerical analysis of \eqref{eq8.3} shows that the wave in its
course becomes steeper in front and less steep behind when $p
<\sqrt{\frac {\bar{h}}{4\sigma}}$  and conversely when $p >
\sqrt{\frac {\bar{h}}{4\sigma}}$.

This result is in contradiction with the assertion of among others
Airy, that a progressive wave always gets steeper in front and
less steep behind. This opinion is conceivable if the dispersion
is neglected.

\subsection{Stability in the Theory of Boussinesq}

Boussinesq considers waves, not necessarily steady, with the same
energy
\begin{equation}\label{eq8.4}
\rho g E = \tfrac 12 \rho g \int_{-\infty}^\infty h^2 dx + \frac
\rho 2 \int_{-\infty}^\infty dx \int_0^{H+h}( u^2+v^2) dy = \rho g
\int_{-\infty}^\infty h^2 dx
\end{equation}
see (\cite{ref9}, pp 85-86).

Furthermore, he introduces the functional
\begin{equation}\label{eq8.5}
  M=\int_{-\infty}^\infty\left\{\left(\frac{\partial h}{\partial
  x}\right)^2 - \frac 3{H ^3} h^3\right\} dx.
\end{equation}
which he calls the ``moment de stabilit\'e'' and he shows that $M$
is a conserved quantity, i.e. independent on $t$, (\cite{ref9}, pp
87-88, 97-99).

After the transformation
$$\eps = \int_x^\infty h^2 dx$$
the expression for $M$ becomes
\begin{equation}\label{eq8.6}
M=\int_0^E \left\{ \left( \frac 14 \frac{\partial h^2}{\partial
\eps}\right)^2 - 3 \frac h{H^3} \right\} d\eps.
\end{equation}

Boussinesq uses, without reference to Euler-Lagrange, the
well-known method to obtain a condition for $h(\eps,t)$ in order
that $M$ attains an extremal value and the result is
\begin{equation}\label{eq8.7}
  1+\frac{2 H^3}{3} h\frac{\partial }{\partial
  \eps}\left( h \frac{\partial h}{\partial
  \eps}\right) = 0
\end{equation}

{}From $d\eps = - h^2 d x = h d\sigma$  follows equation
\eqref{eq4.10} with $\frac{d h}{dt}=0$ and therefore only the
stationary solitary wave with given energy $E$ yields an extremum
for $M$. Variation of $h$ with $\Delta h$ gives $\Delta M> 0$
for all $h(\eps,t)$ and so the extremum of $M$ is a minimum. The
stability of the wave is evident, because also $\Delta M$ does not
depend on $t$.

As in discrete mechanical systems conserved quantities are of
fundamental importance, also in continuous dynamical systems.
Besides the integral invariants $Q = \int_{-\infty}^\infty h \,
dx$ and $E = \int_{-\infty}^\infty h^2 \, dx$, corresponding with
the conservation of mass and energy, Boussinesq discovered a third
invariant, the ``moment de stabilit\'e''.

He also showed that the velocity of the centre of gravity of a
wave does not depend on time and this implies a fourth invariant
of the motion (\cite{ref9}, pp 83-84; \cite{ref16}, p 135).

Conserved functionals may be considered as Hamiltonians in
continuous dynamical systems and they play there a role analogous
to the Hamilton functions in discrete systems. These continuous
dynamical systems have been investigated only rather recently; the
first fundamental results have been established in the seventies
by Lax \cite{ref27}, Gardner \cite{ref28}, Zacharov \cite{ref29}
and Broer \cite{ref30}. Nowadays, there exists an extensive
literature on this subject; a valuable introduction with many
references is the textbook by P.J.~Olver \cite{ref31}.

The Korteweg-de Vries equation is the prototype of an integrable
system with an infinite number of conserved functionals, mutually
in involution with respect to a suitably defined Poisson bracket.
In particular the Korteweg-de Vries equation may be represented
as a Hamiltonian system in the form
\begin{equation}\label{eq8.8}
  \frac{\partial h}{\partial
  t}= - \sqrt{gH} \frac{\partial }{\partial
  x}\delta_h ({\mathcal{H}})
\end{equation}
with $$ {\mathcal{H}}(h) = \int_{-\infty}^\infty \left[ \frac 12
h^2 + \eps \left\{ \left(\frac{\partial h}{\partial x}\right)^2 -
\frac {3}{H^3} h^3 \right\} \right] dx; $$

$\delta_h{\mathcal{H}}$ is the variational derivative of the
Hamilton functional ${\mathcal{H}}$ and $\eps$ a scale parameter,
ref. [13 part I, ch2; part II, ch 5].

The first term in ${\mathcal{H}}$ is the Hamilton functional for
waves in the Lagrange approximation and the second term is the
Boussinesq correction, given by the ``moment de stabilit\'e''
$M$.

Hamilton's theory for finite discrete systems dates from about
1835 and it was a century after Boussinesq that this theory has
been generalized for continuous systems. Boussinesq has set, by
using functionals, a first step into the direction of this
generalization.

\section{Concluding Remarks}

We have discussed in the preceding sections the more important
aspects of the work of Boussinesq and Korteweg-de Vries, who have
besides these also studied other specific topics such as the
velocity field, the path of the fluid particles and the motion of
the centre of gravity of a solitary wave.

Boussinesq finishes his article in the Journal de Math\'ematiques
Pures et Appliqu\'ees with a qualitative examination of the change
of form of long non-stationary waves and an attempt to prove that
a positive solitary wave can be splitted into several other
solitary waves. Korteweg and De Vries wanted to show that their
approximation of the surface of a steady wave may be improved
indefinitely, resulting in a convergent series. They claim in the
introduction of their paper that {\it``in a frictionless liquid
there may exist absolutely stationary waves and that the form of
their surface and the motion of the liquid below it may be
expressed by means of rapidly convergent series''}. The
calculations, however elementary, are so complicated and tedious
that one may expect that these have not received much attention.
Even the second approximation (p. 443), following on the formulae
\eqref{eq6.5} and \eqref{eq7.11} of the present paper requires
already so much effort that it is reasonable to be content with
the first approximation as given in \eqref{eq6.5} and
\eqref{eq7.11}.

It is somewhat surprising that Korteweg and De Vries refer in
their paper only to Boussinesq's short communication in the
Comptes Rendus of 1871 \cite{ref7} and not to the extensive
article in the J. Math. Pures et Appl. \cite{ref9} and the ``Essai
sur la th\'eorie des eaux courantes'' \cite{ref10} in 1872,
respectively 1877. However, we should realize that the
international exchange of scientific achievements in those days
was not at the level as it is today.

As to the credit of the ``a priori demonstration a posteriori'' of
the stable solitary wave, this credit belongs, of course, to M.
Boussinesq. On the other hand, Korteweg and De Vries merit to be
acknowledged for removing doubts on the existence of the ``Great
Wave'' and for their contribution to the theory of long waves in
shallow water.

\subsection*{Acknowledgement}

The author is indebted to the grandsons of Gustav de Vries for
presenting him with a copy of the doctoral thesis of their
grandfather and for records of the handwritten correspondence
between Korteweg and De Vries. He thanks Dr. B.~Willink, historian
at the Erasmus University of Rotterdam and a relative of Korteweg,
for many personal data of Korteweg and De Vries, and for sending
literature, relevant to the content of this essay.

\end{document}